\documentclass{amsart}

\usepackage{amsthm}
\usepackage{amssymb}
\usepackage{amsmath} 

\newcommand{\bea}{\begin{eqnarray*}}
\newcommand{\eea}{\end{eqnarray*}}
\newcommand{\bean}{\begin{eqnarray}}
\newcommand{\eean}{\end{eqnarray}}
\newtheorem{thm}{Theorem}[section]
\newtheorem{prop}[thm]{Proposition}
\newtheorem{lem}[thm]{Lemma}
\newtheorem{cor}[thm]{Corollary}
\newtheorem{rem}[thm]{Remark}
\newtheorem{defn}[thm]{Definition}
\newtheorem{conj}[thm]{Conjecture}
\newtheorem{ex}[thm]{Example}

\newcommand{\HH}{\mathcal{H}_n}

\begin{document}


\title[Subrepresentations in the Polynomial Representation]
{Subrepresentations in the Polynomial Representation
 of the Double Affine Hecke Algebra of type $GL_n$ at $t^{k+1}q^{r-1}=1$}

\author{Masahiro Kasatani}
\email{kasatani@math.kyoto-u.ac.jp}
\address{Department of Mathematics, Graduate School of Science, 
Kyoto University, Kyoto 606-8502, Japan}

\begin{abstract}
We study a Laurent polynomial representation $V$ of
 the double affine Hecke algebra of type $GL_n$ for specialized parameters
 $t^{k+1}q^{r-1}=1$.
We define a series of subrepresentations of $V$
 by using a vanishing condition.
For some cases, we give an explicit basis of the subrepresentation
 in terms of nonsymmetric Macdonald polynomials.
These results are nonsymmetric versions of \cite{FJMM} and \cite{KMSV}.
\end{abstract}

\maketitle

\section{Introduction}

In 1990's, Cherednik introduced the double affine Hecke algebra \cite{C92}.
In terms of the polynomial representation of the algebra,
 Macdonald's conjecture was solved in \cite{C95_Norm}, \cite{C95_DualEval}
 for reduced root systems
 (and in \cite{Sa},\cite{St} for non-reduced $(C^\vee,C)$ case).
For the root system of type $A$,
 a classification of irreducible representations of a certain class
 is given in \cite{C01}, \cite{Su}, \cite{V}.
In \cite{C04}, it is shown that
 finite-dimensional quotients of the polynomial representation
 by the kernel of some degenerate bilinear form are irreducible.

The double affine Hecke algebra $\HH$ of type $GL_n$ is an associative algebra
 with two parameters $t,q$ generated by $X_i$, $Y_i$ and $T_j$
 ($1\leq i \leq n$, $1\leq j \leq n-1$) with some relations.
The algebra $\HH$ has a basic representation $U$
 on the ring of $n$-variable Laurent polynomials.
The representation $U$ is irreducible and {\it $Y$-semisimple}, namely
 the operators $Y_i$ are simultaneously diagonalizable on $U$.
The nonsymmetric Macdonald polynomial is defined
 to be a monic simultaneous eigenvector for $Y_i$.

In this paper, we specialize the parameters at $t^{k+1}q^{r-1}=1$
 for $1\leq k\leq n-1$ and $2\leq r$.
To be precise, introduce a new parameter $u$ and specialize
\bean
(t,q)=(u^{(r-1)/M}, \tau u^{-(k+1)/M}).
 \label{eq:Param}
\eean
Here, $M$ is the greatest common divisor of ($k+1,r-1$)
 and $\tau=\exp(\frac{2\pi\sqrt{-1}}{r-1})$.
We denote by $\HH^{(k,r)}$ and $V$,
the corresponding algebra and its polynomial representation.
The representation $V$ can have subrepresentations
 and they may not be $Y$-semisimple.

In \cite{FJMM} and \cite{KMSV},
 a series of ideals in the ring of symmetric polynomials with $n$-variables
 are defined by vanishing conditions,
 and explicit bases of the ideals are given in terms of
 symmetric Macdonald polynomials specialized at (\ref{eq:Param}).
The vanishing conditions for symmetric polynomial $f$ are
 as follows.
Fix $m\geq1$.
\bea
&&f(x_1,\cdots,x_n)=0 \quad\mbox{if $x_{i_{l,a+1}}=x_{i_{l,a}}tq^{s_{l,a}}$ ($1\leq l \leq m$,$1\leq a\leq k$)}\\
&&\mbox{where $i_{l,a}$ are distinct, $s_{j,i}\in \mathbb{Z}_{\geq0}$, $\sum_{i=1}^{k} s_{j,i}\leq r-1$}.
\eea
This is called the {\it wheel condition} for symmetric case.
In the case $m=1$,
 the basis of the ideal is given in \cite{FJMM}
 by Macdonald polynomials $P_\lambda$ specialized at (\ref{eq:Param})
 with partitions $\lambda$ satisfying $\lambda_i-\lambda_{i+k}\geq r$
 for $1\leq i \leq n-k$.
In the case $(k,r,m)=(1,2,2)$, the basis of the ideal is given in \cite{KMSV}
 by linear combinations of Macdonald polynomials at (\ref{eq:Param}).
Also, a similar ideal in the ring of $BC$-symmetric Laurent polynomials
 is investigated in \cite{K}.

These ideals are invariant under the multiplication
 by symmetric polynomials and Macdonald's $q$-difference operators.
The former actions are symmetric polynomials of $X_i$ and
the latter actions are symmetric polynomials of $Y_i$.
Moreover, the action of $T_i$ on any symmetric polynomial
 is a multiplication by a scalar.
Hence the ideals are representations of the subalgebra of $\HH^{(k,r)}$
 generated by $\{ \mathbb{C}(u)[X_1,\cdots,X_n]^{\mathfrak{S}_n}$,
 $\mathbb{C}(u)[Y_1,\cdots,Y_n]^{\mathfrak{S}_n}$,
 $T_1,\cdots,T_{n-1} \}$.

In this paper, we consider a nonsymmetric version of these ideals.
In other words, we construct a finite series of subrepresentations in $V$
 of the whole algebra $\HH^{(k,r)}$.
In order to obtain them, we define a vanishing condition as follows.
Fix $m\geq1$.
\bean
&&\begin{array}{l}
f(x_1,\cdots,x_n)=0 \quad\mbox{if $x_{i_{l,a+1}}=x_{i_{l,a}}tq^{s_{l,a}}$ ($1\leq l \leq m$,$1\leq a\leq k$)}\\
\mbox{where $i_{l,a}$ are distinct, $s_{l,a}\in \mathbb{Z}_{\geq0}$,}\\
\mbox{ $\sum_{i=1}^{k} s_{l,a}\leq r-2$, and $i_{l,a}<i_{l,a+1}$ if $s_{l,a}=0$}.\label{eq:intro_multiwheel}
\end{array}
\eean
We call the vanishing condition (\ref{eq:intro_multiwheel})
 the {\it wheel condition} for nonsymmetric case.
Denote by $I_m^{(k,r)}$ the space of Laurent polynomials
 satisfying the wheel condition (\ref{eq:intro_multiwheel}).

Let us state the main theorem in this paper.
Define $B^{(k,r)}=\{\lambda\in\mathbb{Z}^n;$
for any $1\leq a \leq n-k$,
$\lambda_{i_{a}}-\lambda_{i_{a+k}}\geq r,$ or
$\lambda_{i_{a}}-\lambda_{i_{a+k}}= r-1$ and $i_{a}<i_{a+k}$ $\}$.
Here, we determine the index $(i_1,\cdots,i_n)=w\cdot(1,\cdots,n)$
 using the shortest element $w\in W=\mathfrak{S}_n$
 such that $\lambda=w\cdot\lambda^+$
 where $\lambda^+$ is the dominant element in the orbit $W\lambda$.
The result is
\begin{thm}\label{thm:main_intro}
The ideal $I_1^{(k,r)}$ is an irreducible representation of $\HH^{(k,r)}$
 and it is $Y$-semisimple.
For any $\lambda\in B^{(k,r)}$,
the nonsymmetric Macdonald polynomial $E_\lambda$ has no pole at $(\ref{eq:Param})$.
Moreover, a basis of the ideal $I_1^{(k,r)}$ is given by
 $\{E_\lambda; \lambda\in B^{(k,r)} \}$ specialized at $(\ref{eq:Param})$.
\end{thm}

For the proof of Theorem \ref{thm:main_intro}, we first show that
these polynomials have no pole at (\ref{eq:Param})
 and they satisfy the wheel condition if they are specialized at (\ref{eq:Param}).
We use the duality relation for nonsymmetric Macdonald polynomials
 (see (\ref{eq:duality}) in Proposition \ref{prop:duality})
 and we count the order of poles and zeros in order to check the statement.
This gives a lower estimate of the character of the ideal.

Next, we give an upper estimate of the character of the ideal.
We introduce the filtration
$V_{(M)}=\mathrm{span}\{x^\lambda ; \lambda\in\mathbb{Z}_n,|\lambda_i|\leq M \}$
 and define a non-degenerate pairing between $V_{(M)}$ and the $n$-th tensor space
 $R_{M,n}=(\mathrm{span}\{e_d; -M\leq d \leq M\})^{\otimes n}$.
We give a spanning set of the quotient space of $R_{M,n}$
 which has the same character as $I^{(k,r)}_1\cap V_{(M)}$.

Finally, we show irreducibility by using intertwining operators.
These operators send one eigenvector to another eigenvector.
We show that $E_\lambda$ for any $\lambda\in B^{(k,r)}$
 is a cyclic vector of $I_1^{(k,r)}$.

For the case $n=k+1$, we also show that $V/I^{(n-1,r)}_1$ is irreducible
 and we give an explicit basis of $V/I^{(n-1,r)}_1$.
We expect that all subquotients $I^{(k,r)}_m/I^{(k,r)}_{m-1}$
 of the series are irreducible.

The plan of the paper is as follows.
In Section 2, we review the double affine Hecke algebra,
 the polynomial representation $U$,
 the nonsymmetric Macdonald polynomials $E_\lambda$,
 and intertwiners.
In Section 3, we state the wheel condition and
 show that it determines a subrepresentaion $I_1^{(k,r)}$.
In Section 4,
 we give a lower estimate of the character of $I_1^{(k,r)}$
 using nonsymmetric Macdonald polynomials.
In Section 5, we give an upper estimate of the character of $I_1^{(k,r)}$,
 and we show that $I_1^{(k,r)}$ is irreducible.
In Section 6, we define a (finite) series of subrepresentations $I_m^{(k,r)}$
 by the wheel condition.
The series contains the irreducible representation $I_1^{(k,r)}$
 defined in Section 3.
We also treat the case $n=k+1$ in Section 6.

\subsection*{Acknowledgments.}
The author thanks Ivan Cherednik for a useful discussion
 which leads him to this study.
The author also thanks his adviser Tetsuji Miwa
 for giving many comments to the manuscript.

\section{Double affine Hecke algebra and nonsymmetric Macdonald polynomials}

In this section, we review the double affine Hecke algebra of type $GL_n$,
 nonsymmetric Macdonald polynomials, and intertwiners.
For more details, see, e.g., \cite{C95_nonsym}.
In this paper, we follow the notation in \cite{MN}.

\subsection{Double affine Hecke algebra of type $GL_n$}

Let $\tilde{\mathbb{K}}$ be the field $\mathbb{C}(t^{1/2},q)$.

\begin{defn}\normalfont
The double affine Hecke algebra $\HH$ of type $GL_n$ is
 an associative $\tilde{\mathbb{K}}$-algebra generated by
\[
\langle X_1^{\pm1},\cdots,X_n^{\pm1},Y_1^{\pm1},\cdots,Y_n^{\pm1},
 T_1,\cdots,T_{n-1} \rangle
\]
satisfying the following relations:
\begin{eqnarray*}
&& \mbox{$X_i^{\pm1}$ are mutually commutative and $X_iX_i^{-1}=1$}, \\
&& \mbox{$Y_i^{\pm1}$ are mutually commutative and $Y_iY_i^{-1}=1$}, \\
&& (T_i-t^{1/2})(T_i+t^{-1/2})=0, \\
&& T_iT_{i+1}T_i=T_{i+1}T_iT_{i+1}, \\
&& T_iT_j=T_jT_i \quad(|i-j|\geq2),\\
&& T_iX_iT_i=X_{i+1}, \\
&& T_iX_j=X_jT_i \quad(j\neq i,i+1),\\
&& T_i^{-1}Y_iT_i^{-1}=Y_{i+1}, \\
&& T_iY_j=Y_jT_i \quad(j\neq i,i+1),\\
&& Y_2^{-1}X_1Y_2X_1^{-1}=T_1^2, \\
&& Y_i\tilde{X}=q\tilde{X}Y_i
 \quad\mbox{where $\tilde{X}=\prod_{i=1}^{n} X_i$},\\
&& X_i\tilde{Y}=q^{-1}\tilde{Y}X_i
 \quad\mbox{where $\tilde{Y}=\prod_{i=1}^{n} Y_i$}.
\end{eqnarray*}

\end{defn}

\subsection{Polynomial representation}
The algebra $\HH$ has a basic representation on the ring of Laurent polynomials
 $U=\tilde{\mathbb{K}}[x_1^{\pm1},\cdots,x_n^{\pm1}]$.
\bea
X_i&\mapsto&x_i,  \\
T_i&\mapsto&t^{1/2}s_i+\frac{t^{1/2}-t^{-1/2}}{x_i/x_{i+1}-1}(s_i-1),  \\
Y_i&\mapsto&T_iT_{i+1}\cdots T_{n-1}\omega T_1^{-1}T_2^{-1}\cdots T_{i-1}^{-1}.
\eea
Here, $s_i$ is the permutation of the variables $x_i$ and $x_{i+1}$,
 $\omega=s_{n-1}\cdots s_1\tau_{1}$, and $\tau_ix_j=q^{\delta_{ij}}x_j$.
Namely, $\omega f(x_1,x_2,\cdots,x_n)=f(qx_n, x_1,x_2,\cdots,x_{n-1})$
 for any $f\in U$.

\subsection{Nonsymmetric Macdonald polynomials $E_\lambda$}

The representation $U$ is semisimple with respect to the action of $\{Y_i\}$
and it has simultaneous eigenvectors labeled by $\lambda\in \mathbb{Z}^n$.
We call such eigenvectors {\it $Y$-eigenvectors}.
The monic $Y$-eigenvecrtors are called nonsymmetric Macdonald polynomials.
To describe eigenvalues, we use the following notations.

Denote the standard basis of $\mathbb{Z}^n$ by $\{\epsilon_i\}$.
We identify the weight lattice $P$ of type $GL_n$ with $\mathbb{Z}^n$.
For $\lambda\in P$, let us denote by $\lambda_i$ the $i$-th component of $\lambda$.
Let $\Delta_+=\{\epsilon_i-\epsilon_{j};i<j\}$ be the set of positive roots
 for $A_{n-1}$ and $W=\mathfrak{S}_n$ the Weyl group.
Set $P^+=\{\lambda\in P;\langle \lambda,\alpha_i\rangle\geq 0$ for any $i$$\}$
 and write the dominant element $\lambda^+\in P^+ \cap W\lambda$.
\begin{eqnarray*}
&&\rho=\frac{1}{2}\sum_{\alpha \in \Delta_+}\alpha
=\left(\frac{n-1}{2},\frac{n-3}{2},\cdots,-\frac{n-1}{2}\right),  \\
&&\rho(\lambda)=\frac{1}{2}\sum_{\alpha \in \Delta_+}\chi(\langle \lambda,\alpha \rangle)\alpha,\\
&&\chi(a)=\left\{
\begin{array}{ll}
 +1 & (a\geq0) \\
 -1 & (a<0) \\
\end{array}\right..
\end{eqnarray*}
The element $\rho(\lambda)$ is equal to $w_\lambda\rho$ where $w_\lambda\in W$
 is the shortest element such that $\lambda=w_\lambda\lambda^+$.
Note that $\rho(\lambda)_i>\rho(\lambda)_j$ if ($\lambda_i>\lambda_j$) or
($\lambda_i=\lambda_j$ and $i<j$).

We define the ordering $\succ$ by $\lambda \succ \mu \Leftrightarrow
 (\lambda^+>\mu^+)$ or $(\lambda^+=\mu^+$ and $\lambda>\mu)$.
Here $>$ is the dominance ordering: $\lambda\geq\mu \Leftrightarrow 
\sum_{j=1}^{l}\lambda_j \geq\sum_{j=1}^{l} \mu_j$ for any $1\leq l\leq n$.
The operator $Y_i$ is triangular with respect to the ordering $\succ$:
\[Y_ix^\lambda=t^{\rho(\lambda)_i}q^{\lambda_i}x^\lambda
 + \sum_{\mu\prec\lambda} c_{\lambda,\mu}x^\mu.\]

We write $f(t^{\pm\rho(\lambda)}q^{\pm\lambda})=f(t^{\pm\rho(\lambda)_1}q^{\pm\lambda_1},\cdots, t^{\pm\rho(\lambda)_n}q^{\pm\lambda_n})$.
The nonsymmetric Macdonald polynomial $E_\lambda$ is defined to be the monic $Y$-eigenvector:
\bea
f(Y)E_\lambda&=&f(t^{\rho(\lambda)}q^\lambda)E_\lambda \quad\mbox{for any $f\in U$}, \\
E_\lambda&=&x^\lambda+\sum_{\mu\prec\lambda} c_{\lambda,\mu}x^\mu.
\eea
The nonsymmetric Macdonald polynomials form a basis of $U$.

\subsection{Duality relation}
For any $\lambda\in P$ and $f\in U$,
we use the notation
\[ u_\lambda(f)=f(t^{-\rho(\lambda)}q^{-\lambda}). \]
For example, $u_\lambda(x_i)=t^{-\rho(\lambda)_i}q^{-\lambda_i}$.
\begin{prop}[duality\cite{C95_nonsym}]\label{prop:duality}
Two nonsymmetric Macdonald polynomials $E_\lambda$ and $E_\mu$ satisfy
the following relation.
\bean
 \frac{u_\mu(E_\lambda)}{u_0(E_\lambda)}
=\frac{u_\lambda(E_\mu)}{u_0(E_\mu)}. \label{eq:duality}
\eean
\end{prop}

\subsection{Recurrence formulas for $u_0(E_{\lambda})$}
We have the following recurrence formulas for
 $u_0(E_\lambda)=E_{\lambda}(t^{-\rho})$.
\begin{lem}\label{lem:u_0formula_omega}
We write $\omega\lambda=(\lambda_2,\cdots,\lambda_n,\lambda_1+1)$.
Then we have
\[
u_0(E_{\omega\lambda})=u_0(E_{\lambda})\cdot t^{\rho(\lambda)_1}.
\]
\end{lem}
\begin{lem}\label{lem:u_0formula_s_i}
Suppose $\lambda \prec s_i\lambda \ (namely, \langle\lambda,\alpha_i\rangle<0)$.
Then we have
\bea
u_0(E_{s_i\lambda})=u_0(E_{\lambda})\cdot
t^{-1}\frac{tu_\lambda(x_i/x_{i+1})-1}{u_\lambda(x_i/x_{i+1})-1}.
\eea
\end{lem}

\subsection{Intertwiners}
There exist some operators which create $E_{\omega\lambda}$ and $E_{s_i\lambda}$ from $E_\lambda$.
In other words, they intertwine the eigenvectors.
This is due to Cherednik (see e.g. \cite{C95_nonsym}).
Especially for $GL_n$ case, see also \cite{Kn} and \cite{MN}.

\begin{lem}[operator $A$]
We set $\omega=(T_1\cdots T_{n-1})^{-1}Y_1$ and $A=\omega X_1$.
Then we have,
\bea
A E_\lambda &=& q^{\lambda_1+1} E_{\omega\lambda}.
\eea
\end{lem}

\begin{lem}[operator $B_i$]\label{lem:BActionForGeneric}
We set
\bea
B_i(\lambda)=T_i+\frac{t^{1/2}-t^{-1/2}}{u_\lambda(x_i/x_{i+1})-1}.
\eea
Note that for any $\lambda\in P$, the denominator is not zero.
For simplicity, we write $B_i$ when $B_i(\lambda)$ acts on $E_\lambda$.
\bea
\mbox{If }s_i \lambda \succ \lambda, &\mbox{then}& B_i E_\lambda = t^{1/2}E_{s_i\lambda}.\\
\mbox{If }s_i \lambda = \lambda, &\mbox{then}& B_i E_\lambda = 0. \\
\mbox{If }s_i \lambda \prec \lambda, &\mbox{then}& B_i E_\lambda = t^{-1/2}\\
&&\quad \times
 \frac{(t^{-1}u_\lambda(x_i/x_{i+1})-1)
 (tu_\lambda(x_i/x_{i+1})-1)}
 {(u_\lambda(x_i/x_{i+1})-1)^2}
E_{s_i\lambda}.
\eea
\end{lem}

Note that the group action generated by $\omega^{\pm1}$ and $s_i$
 on $P$ is transitive.
Hence by applying the intertwiners $A$ and $B_i$,
 we see that $E_\lambda$ is a cyclic vector in $U$ for any $\lambda\in P$.
Namely, $U$ is irreducible.

\section{Irreducible representation defined by wheel condition}

In this section, we impose the specialization of parameters (\ref{eq:Param}).
Namely, let $k,r$ be integers with $n-1\geq k\geq 1$ and $r\geq 2$,
and we specialize parameters $t^{k+1}q^{r-1}=1$.
To be precise, introduce a new parameter $u$ and specialize
$(t,q)=(u^{(r-1)/M}, \tau u^{-(k+1)/M})$.
Here, $M$ is the greatest common divisor of ($k+1$, $r-1$)
 and $\tau=\exp(\frac{2\pi\sqrt{-1}}{r-1})$.

Take
\bea
\mathbb{K}'&=&\{ c\in\tilde{\mathbb{K}} ; \mbox{ $c$ is regular at (\ref{eq:Param}) } \},\\
{\HH}'&=&\{ h\in\HH ; \mbox{ $h$ is regular at (\ref{eq:Param}) } \},\\
U'&=&\{ f\in U ; \mbox{ $f$ is regular at (\ref{eq:Param}) } \},
\eea
and let $\mathbb{K}$, $\HH^{(k,r)}$ and $V$ be
 the images of $\mathbb{K}'$, ${\HH}'$, and $U'$
 by the specialization (\ref{eq:Param}).
Note that $\mathbb{K}=\mathbb{C}(u)$
and $V=\mathbb{K}[x_1^{\pm1},\cdots,x_n^{\pm1}]$.

In this situation,
 we construct an ideal of $V$ defined by a certain vanishing condition.
This gives an irredicble representation of $\HH^{(k,r)}$.

\begin{defn}\normalfont\label{defn:WheelCondition}
Define
\bea
Z^{(k,r)}=\{ (z_1,\cdots,z_n)\in\mathbb{K}^n&;&
\mbox{there exist distinct $i_1,\cdots,i_{k+1}\in\{1,\cdots,n\}$}\\
&&\mbox{and positive integers $s_1,\cdots, s_{k+1}\in \mathbb{Z}_{\geq0}$}\\
&&\mbox{ such that $z_{i_{a+1}}=z_{i_a}tq^{s_a}$ for $1\leq a\leq k $,}\\
&&\mbox{$\sum_{a=1}^{k} s_a\leq r-2$, and $i_a<i_{a+1}$ if $s_a=0 \}$.}
\eea
We define the ideal
\[
I^{(k,r)}=\{f\in V ; f(z)=0 \mbox{ for any $z\in Z^{(k,r)}$} \}.
\]
We call the defining condition of $I^{(k,r)}$ the {\it wheel condition}.
\end{defn}

\begin{rem}\normalfont
Let us denote the relation $z_j=tq^{s}z_i$ by $z_i \stackrel{s}{\rightarrow} z_j$.
Take an element $z\in Z^{(k,r)}$.
Then there exist $({i_1},\cdots,{i_{k+1}})$
 such that\\
\begin{tabular}{rcl}
$z_{i_{1}}$ & $\stackrel{s_1}{\rightarrow} z_{i_{2}} \stackrel{s_2}{\rightarrow}$ & $\cdots$ \\
${}_{s_{k+1}} \uparrow$ &  & $\downarrow {}_{s_{k-1}}$ \\
$z_{i_{k+1}}$ & $\stackrel{s_k}{\leftarrow}$ & $z_{i_{k}}$
\end{tabular}.\\
Note that $s_{k+1}:=r-1-\sum_{a=1}^{k} s_a\geq 1$,
 and the specialization $t^{k+1}q^{r-1}=1$ implies that
 $z_{i_{k+1}} \stackrel{s_{k+1}}{\rightarrow} z_{i_1}$.
It looks like a wheel.
This is the reason why we call such a condition the {\it wheel} condition.
\end{rem}

\begin{rem}\normalfont
The wheel condition is originally appeared in \cite{FJMM}
 as a vanishing condition for symmetric polynomials
on the set
$Z_{sym}^{(k,r)}=\{ (z_1,\cdots,z_n)\in\mathbb{K}^n;$
$\exists s_1,\cdots, s_{k}\in \mathbb{Z}_{\geq0}$ such that
$z_{a+1}=z_{a}tq^{s_a}$ for $1\leq a\leq k $,
$\sum_{a=1}^{k} s_a\leq r-1 \}$.
The set $Z_{sym}^{(k,r)}$ is apparently different from $Z^{(k,r)}$.
However, for symemtric polynomials,
 the vanishing condition on $Z_{sym}^{(k,r)}$ is equivalent to that on $Z^{(k,r)}$.
\end{rem}

Now we describe the first main statement.

\begin{prop}\label{prop:rep}
The ideal $I^{(k,r)}$ is a representation of the algebra $\HH^{(k,r)}$.
\end{prop}

Before giving the proof,
let us give equivalent definitions of the ideal $I^{(k,r)}$.
We can reduce the set $Z^{(k,r)}$ to smaller subsets.

\begin{defn}\normalfont
Let $\lambda\in P$.
Suppose that $(u_\lambda(x_1),\cdots,u_\lambda(x_n))\in Z^{(k,r)}$.
Then there exist $(i_1,\cdots,i_{k+1})$
 and $s_1,\cdots s_{k} \in \mathbb{Z}_{\geq 0}$ satisfying
\bea
\mbox{$u_\lambda(x_{i_{a+1}})=u_\lambda(x_{i_{a}})tq^{s_a}$ for $1\leq a\leq k$},\\
\mbox{$\sum_a s_a\leq r-2$, and $i_a<i_{a+1}$ if $s_a=0$}.
\eea
We call such $(i_1,\cdots,i_{k+1})$ {\it a wheel in $\lambda$}.
For some $\sigma\in \mathfrak{S}_{k+1}$,
if $(i_1,\cdots,i_{k+1})$ and
 $\sigma(i_1,\cdots,i_{k+1})=(i_{\sigma^{-1}(1)},\cdots,i_{\sigma^{-1}(k+1)})$
 are wheels in $\lambda$, we identify them.
In such a case, we see that $\sigma(i_1,\cdots,i_{k+1})=$
$(i_{a},i_{a+1},\cdots,i_{k+1},i_1,\cdots,{i_{a-1}})$
 for some $1 \leq a\leq k+1$ by the definition of wheels.
Denote the number of equivalent classes of wheels by $\sharp^{(k,r)}(\lambda)$.
\end{defn}

We shall introduce two subsets $S^{(k,r)}$ and $S'{}^{(k,r)}$ in $P$.

\begin{defn}\normalfont
Let $a,b$ be integers with $a\geq 2$ and $b\geq 1$.
Take $\lambda\in P$.
We call $(i,j)$ is a {\it neighborhood of type $(a,b)$ in $\lambda$} if
(i) and (ii) hold:\\
\ (i) $\rho(\lambda)_i-\rho(\lambda)_j=a-1$,\\
\ (ii) ($\lambda_i-\lambda_j\leq b-1$), or ($\lambda_i-\lambda_j= b$ and $j<i$).

\end{defn}

\begin{defn}\normalfont
We define $S^{(k,r)}=\{ \lambda\in P ;$
$\lambda$ has a neighborhood of type $(k+1,r-1)$ $\}$ and
 $S'{}^{(k,r)}=\{ \lambda\in S^{(k,r)} ;$
$\lambda$ has a neighborhood $(i,j)$ of type $(k+1,r-1)$
 such that $\lambda_i-\lambda_j\leq r-2$ $\}$.
\end{defn}

Fix $\lambda\in S^{(k,r)}$ and
let $(i_1,i_2,\cdots,i_n)=w_\lambda\cdot(1,2,\cdots,n)$.

Suppose that $(i_l,i_{l+k})$ is a neighborhood of type $(k+1,r-1)$ in $\lambda$.
If $\lambda_{i_l}-\lambda_{i_{l+k}}\leq r-2$, then $(i_l,i_{l+1},\cdots,i_{l+k})$
 is a wheel in $\lambda$.
If $\lambda_{i_l}-\lambda_{i_{l+k}}=r-1$ and $i_{l+k}<i_l$, then
 $(i_{l+a},\cdots,i_{l+k},i_l,\cdots,{i_{l+a-1}})$ is a wheel in $\lambda$.
Here, we take $a\geq1$ satisfying $\lambda_{i_l}=\lambda_{i_{l+1}}=\cdots=\lambda_{i_{l+a-1}}>\lambda_{i_{l+a}}$.

Suppose that $(i_l,i_{l+k})$ and $(i_m,i_{m+k})$ are
 different neighborhoods of type $(k+1,r-1)$ in $\lambda$
(namely $l\neq m$).
Consider the wheels in $\lambda$ defined in the previous paragraph.
Then we see that these wheels belong in different equivalent classes of wheels.
Hence $\sharp^{(k,r)}(\lambda)$ is greater than or equal to
 the number of neighborhoods of type $(k+1,r-1)$.

\begin{ex}\normalfont
($n=3,k=1,r=2$): \quad 
The set $S^{(1,2)}$ and  $S'{}^{(1,2)}$ are given as follows:
\bea
S^{(1,2)}=&&\{ (\lambda_0,\lambda_1,\lambda_1),
 (\lambda_1,\lambda_0,\lambda_1),
 (\lambda_1,\lambda_1,\lambda_0), \\
&&\ (\lambda_0,\lambda_1,\lambda_1+1),
 (\lambda_1,\lambda_0,\lambda_1+1),
 (\lambda_1,\lambda_1+1,\lambda_0), \\
&&\qquad ;\lambda_0,\lambda_1\in\mathbb{Z}\}, \\
S'{}^{(1,2)}=&&\{ (\lambda_0,\lambda_1,\lambda_1),
 (\lambda_1,\lambda_0,\lambda_1),
 (\lambda_1,\lambda_1,\lambda_0), \\
&&\qquad ;\lambda_0,\lambda_1\in\mathbb{Z}\}. 
\eea
\end{ex}

\begin{rem}\label{rem:omega}\normalfont
Denote $\overline{0}=n$ and $\overline{i}=i$ for $1\leq i \leq n-1$.
Then $(i_1,\cdots,i_{k+1})$ is a wheel in $\lambda$ if and only if
 $(\overline{i_a-1},\cdots,\overline{i_{k+1}-1},\overline{i_1-1},\cdots,\overline{i_{a-1}-1})$
 is a wheel in $\omega\lambda$ for some $1\leq a\leq k+1$.
We also see that $(i,j)$ is a neighborhood of type $(a,b)$ in $\lambda$
 if and only if $(\overline{i-1},\overline{j-1})$ is a neighborhood
 of type $(a,b)$ in $\omega\lambda$ for any $a\geq 2,b\geq1$.
\end{rem}

\begin{lem}\label{lem:EquivDef}
The ideal $I^{(k,r)}$ coincides with the following ideals:
\bea
&&J_1=\{f\in V; \mbox{$u_\lambda(f)=0$
 for any $\lambda\in P$ satisfying $\sharp^{(k,r)}(\lambda)\geq 1$} \},\\
&&J_2=\{f\in V; \mbox{$u_\lambda(f)=0$ for any $\lambda\in S^{(k,r)}$} \},\\
&&J_3=\{f\in V; \mbox{$u_\lambda(f)=0$ for any $\lambda\in S'{}^{(k,r)}$} \}.
\eea
\end{lem}

\begin{proof}
We see that $\{\lambda\in P; \sharp^{(k,r)}(\lambda)\geq 1\}
\supset S^{(k,r)} \supset S'{}^{(k,r)}$.
Hence $J_1\subset J_2 \subset J_3$.
If $f\in I^{(k,r)}$, then by the definition of $\sharp^{(k,r)}$,
 we see $f\in J_1$.

Let us show that $J_3\subset I^{(k,r)}$.
Fix an element $f=\sum_{\mu\in P} c_\mu x^\mu \in J_3$.
Denote $\max\{|\mu_i|;c_\mu\neq0,1\leq i \leq n\}$ by $\deg(f)$.
 and take an integer $N > 2(\deg(f)+1)$.
Let $M > 2N+2[\frac{n}{k+1}](r-1)$.

Take $\lambda^{(0)}\in S'{}^{(k,r)}$
 satisfying the following condition:
\bean
&&\begin{array}{l}
\mbox{for some $1\leq l \leq n-k$,}\\
\mbox{$\lambda^{(0)}$ has a neighborhood $(i_l,i_{l+k})$ of type $(k+1,r-1)$}\\
\mbox{such that $\lambda^{(0)}_{i_l}-\lambda^{(0)}_{i_{l+k}}\leq r-2$, and
 $|\lambda^{(0)}_i-\lambda^{(0)}_j|> M$}\\
\mbox{ for any $1\leq i\neq j\leq n$ except for $(i,j)=(i_a,i_b)$ ($l\leq  a,b\leq l+k$)}.
\end{array}\label{eq:lambda^0}
\eean
Here we set $(i_1,\cdots,i_n)=w_{\lambda^{(0)}}\cdot(1,\cdots,n)$.
We see that $\sharp^{(k,r)}(\lambda^{(0)})=1$.
Define a finite set
\bea
S(\lambda^{(0)})=\{ \lambda\in P ; \mbox{$0\leq \lambda_i-\lambda^{(0)}_i \leq N$
 for $i\neq i_l,i_{l+1},\cdots,i_{l+k}$} \}.
\eea
We see $S(\lambda^{(0)})\subset S'{}^{(k,r)}$.
Thus $u_\lambda(f)=0$ for $\lambda\in S(\lambda^{(0)})$.
Because $f$ is a Laurent polynomial of $2\deg (f)\leq N$,
we see that $f(z)=0$ for any $z\in\mathbb{K}^n$ satisfying
$z_{i_1}=u_{\lambda^{(0)}}(x_{i_1})$ and
$z_{i_{a+1}}=tq^{\lambda^{(0)}_{i_a}-\lambda^{(0)}_{i_{a+1}}}z_{i_a}
 (1\leq a\leq k)$.

Choose any $\lambda^{(0)}$ satisfying (\ref{eq:lambda^0}).
Then by a similar argument,
 we see that $f(z)=0$ for any $z\in Z^{(k,r)}$.
Therefore $f\in I^{(k,r)}$.
\end{proof}

\begin{proof}[Proof of Proposition \ref{prop:rep}]
We consider the actions of generators $T_i,X_j^{\pm1}$ and $Y_j^{\pm1}$
 ($1\leq i \leq n-1$, $1\leq j \leq n$).
By the definition, we see $X_i^{\pm1} f=x_i^{\pm1}f \in I^{(k,r)}$
 for any $f\in I^{(k,r)}$.
Since $Y_i^{\pm1}$ is a linear combination of products of
 $\omega^{\pm1}$ and $T_i$,
it is sufficient to show that $\omega^{\pm1} f$ and
 $T_i f\in I^{(k,r)}$ ($1\leq i \leq n-1$) for any $f\in I^{(k,r)}$.

Fix an element $f\in I^{(k,r)}$.
Similarly to Lemma \ref{lem:EquivDef},
take an integer $N > 2(\deg(f)+1)$ and
let $M > 2N+2[\frac{n}{k+1}](r-1)$.
Take $\lambda^{(0)}\in S'{}^{(k,r)}$ satisfying (\ref{eq:lambda^0})
and define a finite set
\bea
S(\lambda^{(0)})=\{ \lambda\in P ; \mbox{$0\leq \lambda_i-\lambda^{(0)}_i \leq N$
 for $i\neq i_l,\cdots,i_{l+k}$} \}.
\eea

Let us show
 $u_\lambda(\omega^{\pm1} f)=0$ 
 for any $\lambda\in S(\lambda^{(0)})$.
By the definition of $u_\lambda$ and $\omega$,
 we see $u_\lambda(\omega^{\pm1} f)=u_{\omega^{\pm1}\lambda}(f)$.
Since $\omega^{\pm1}\lambda$ belongs to $S^{(k,r)}$ for any $\lambda\in S(\lambda^{(0)})$,
 we see $(u_{\omega^{\pm1}\lambda}(x_1),\cdots,u_{\omega^{\pm1}\lambda}(x_n))\in Z^{(k,r)}$.
Hence we have $u_\lambda(\omega^{\pm1} f)=0$.

Let us show
 $u_\lambda(T_i f)=0$ ($1\leq i \leq n-1$)
 for any $\lambda\in S(\lambda^{(0)})$.
Recall $T_i=-\frac{t^{1/2}-t^{-1/2}}{x_i/x_{i+1}-1}+t^{-1/2}\frac{tx_i/x_{i+1}-1}{x_i/x_{i+1}-1}s_i$.
Since $u_\lambda(f)=0$, it is easy to see that
\bea
u_\lambda\left(-\frac{t^{1/2}-t^{-1/2}}{x_i/x_{i+1}-1}f\right)=0.
\eea
Let us show that
\bean
u_\lambda
\left(t^{-1/2}\frac{tx_i/x_{i+1}-1}{x_i/x_{i+1}-1}s_if\right) =0.
 \label{eq:2ndCoeff}
\eean
If $\lambda_i=\lambda_{i+1}$, then (\ref{eq:2ndCoeff}) holds
 because $u_\lambda(tx_i/x_{i+1})=1$.
If $\lambda_i\neq\lambda_{i+1}$, then $s_i\lambda\in S'{}^{(k,r)}$.
Therefore $u_\lambda(s_if)=u_{s_i\lambda}(f)=0$.
Hence (\ref{eq:2ndCoeff}) is proved.

We have proved that
 $u_\lambda(\omega^{\pm1} f)=u_\lambda(T_i f)=0$ ($1\leq i \leq n-1$)
 for $\lambda\in S(\lambda^{(0)})$ and for any choice of
 $\lambda^{(0)}\in S'{}^{(k,r)}$ satisfying (\ref{eq:lambda^0}).
Note that $2\deg(\omega f), 2\deg(T_i f)\leq N$.
Hence from the same argument as Lemma \ref{lem:EquivDef},
 we see that $\omega^{\pm1} f$ and $T_i f\in I^{(k,r)}$ ($1\leq i \leq n-1$).
Therefore the desired statement is proved.
\end{proof}

Now we come to the main theorem of the paper.
Set $B^{(k,r)}=P\backslash S^{(k,r)}$.

\begin{thm}\label{thm:NoWheel}
The ideal $I^{(k,r)}$ is irreducible.
For any $\lambda\in B^{(k,r)}$, the nonsymmetric Macdonald
 polynomial $E_\lambda$ has no pole at the specialization $(\ref{eq:Param})$.
A basis of $I^{(k,r)}$ is given by
 $\{ E_\lambda ; \lambda\in B^{(k,r)} \}$
 specialized at $(\ref{eq:Param})$.
\end{thm}

We give a proof of the theorem in the next section.

\section{Proof of Theorem \ref{thm:NoWheel}}

We use the following notation for multiplicity of zeros and poles.
\begin{defn}\normalfont
Let $M$ be the greatest common divisor of $k+1$ and $r-1$.
For $c\in\tilde{\mathbb{K}}$,
denote by $\zeta(c)$ the integer satisfying
\[
c=\left(t^{(k+1)/M}q^{(r-1)/M}-e^{2\pi\sqrt{-1}/M}\right)^{\zeta(c)}c'
\]
 where $c'\in\tilde{\mathbb{K}}$ have no zero or pole
 at (\ref{eq:Param}).
\end{defn}
Note that $\zeta(c)$ is the order of zeros or poles of $c$
 at the relevant irreducible component of $t^{k+1}q^{r-1}=1$.

To prepare for the proof, we give a key lemma.
This lemma claims that the changes of $\sharp^{(k,r)}(\lambda)$
 and $\zeta(u_0(E_\lambda))$ when $s_i$ acts on $\lambda$ are related to
 the value of $u_\lambda(x_i/x_{i+1})$.

\begin{lem}[Key Lemma]\label{lem:KeyLem1}
Suppose $s_i\lambda\succ\lambda$.\\
$(i)$
If $\sharp^{(k,r)}(s_i\lambda)<\sharp^{(k,r)}(\lambda)$,
 then $u_\lambda(x_i/x_{i+1}) = t^{-1}$. \\
\ \quad If $\sharp^{(k,r)}(s_i\lambda)>\sharp^{(k,r)}(\lambda)$,
 then $u_\lambda(x_i/x_{i+1}) = t$.\\
$(ii)$
We have $|\zeta(u_0(E_{s_i\lambda}))-\zeta(u_0(E_\lambda))|\leq 1$.
Moreover,
\begin{eqnarray*}
 u_\lambda(x_i/x_{i+1}) = t^{-1}
 &\Leftrightarrow& \zeta(u_0(E_{s_i\lambda}))=\zeta(u_0(E_\lambda))+1 ,\\
 u_\lambda(x_i/x_{i+1}) = 1
 &\Leftrightarrow& \zeta(u_0(E_{s_i\lambda}))=\zeta(u_0(E_\lambda))-1.
\end{eqnarray*}
\end{lem}

\begin{rem}\normalfont
In fact, the converse statements for (i) are true.
Here we omit them because we do not use them for the proof of Theorem \ref{thm:NoWheel}.
\end{rem}

\begin{proof}
(i).\quad
Fix $(i_1,\cdots,i_{k+1})\in \{1,\cdots,n\}^{k+1}$ ($i_a$ are distinct).

Suppose $\{i,i+1\}\not\subset\{i_1,\cdots,i_{k+1}\}$. Then,
$(i_1,\cdots,i_{k+1})$ is a wheel in $\lambda$ if and only if
$s_i(i_1,\cdots,i_{k+1})$ is a wheel in $s_i\lambda$.

Suppose $\{i,i+1\}\subset\{i_1,\cdots,i_{k+1}\}$.
(Case a) If $(i_1,\cdots,i_{k+1})$ is a wheel in $\lambda$ and
 $s_i(i_1,\cdots,i_{k+1})$ is not a wheel in $s_i\lambda$,
 then by the definition of the wheel,
 $u_\lambda(x_{i+1})$ must be equal to $u_\lambda(x_{i})tq^0$.
 Thus $u_{i\lambda}(x_i/x_{i+1})=t^{-1}$.
(Case b) If $(i_1,\cdots,i_{k+1})$ is not a wheel in $\lambda$ and
 $s_i(i_1,\cdots,i_{k+1})$ is a wheel in $s_i\lambda$,
 then by the same reason as (a), we see $u_{s_i\lambda}(x_i/x_{i+1})=t^{-1}$.
Note that $u_\lambda(x_i/x_{i+1})=u_{s_i\lambda}(x_i/x_{i+1})^{-1}$.
Hence if the case (a) (resp. (b)) occurs for $(i_1,\cdots,i_{k+1})$,
then the case (a) (resp. (b)) does not occur for any other
 $(i'_1,\cdots,i'_{k+1})$ such that $\{i,i+1\}\subset\{i'_1,\cdots,i'_{k+1}\}$.
Namely two cases (a) and (b) do not concur.

Hence, if $\sharp^{(k,r)}(s_i\lambda)<\sharp^{(k,r)}(\lambda)$,
 then the case (a) does occur and the case (b) does not occur.
Thus $u_\lambda(x_i/x_{i+1})=t^{-1}$.
If $\sharp^{(k,r)}(s_i\lambda)>\sharp^{(k,r)}(\lambda)$,
  then the case (b) does occur and the case (a) does not occur.
Thus $u_\lambda(x_i/x_{i+1})=t$.

(ii).\quad
This is obvious from Lemma \ref{lem:u_0formula_s_i}.
\end{proof}

In the proof above, we see the following fact.
\begin{rem}\normalfont\label{rem:GenAndBrakeOfWheel}
The relation $\sharp^{(k,r)}(\lambda)>\sharp^{(k,r)}(s_i\lambda)$
 holds if and only if
 $(i_1,\cdots,i_{k+1})$ is a wheel in $\lambda$ and
 $s_i(i_1,\cdots,i_{k+1})$ is not a wheel in $s_i\lambda$
 for some $(i_1,\cdots,i_{k+1})$.
Conversely, the relation $\sharp^{(k,r)}(\lambda)<\sharp^{(k,r)}(s_i\lambda)$
 holds if and only if
 $(i_1,\cdots,i_{k+1})$ is not a wheel in $\lambda$ and
 $s_i(i_1,\cdots,i_{k+1})$ is a wheel in $s_i\lambda$
 for some $(i_1,\cdots,i_{k+1})$.
\end{rem}

\begin{defn}\normalfont
We introduce an equivalence relation {\it"intertwined"}
 generated by the following two relations:\\
(i) We call $\lambda$ and $\omega\lambda$ are {\it intertwined}.\\
(ii) We call $\lambda$ and $s_i\lambda$ for $s_i\lambda\neq\lambda$ are
   {\it intertwined} if $u_\lambda(x_i/x_{i+1})\neq1,t^{\pm1}$.
\end{defn}
From Lemma \ref{lem:u_0formula_omega}, Remark \ref{rem:omega},
 and Key Lemma \ref{lem:KeyLem1},
 if $\lambda$ and $\lambda'$ are intertwined, then we see that
 $\zeta(u_0(E_{\lambda}))=\zeta(u_0(E_{\lambda'}))$ and
 $\sharp^{(k,r)}(\lambda)=\sharp^{(k,r)}(\lambda')$.
Moreover, if $E_\lambda$ have no pole at (\ref{eq:Param}),
 then by applying intertwiners $A$ and $B_i$, we obtain $E_{\lambda'}$.

There is a connection between the number of neighborhoods of different types.

\begin{lem}\label{lem:NbdOfDiffTypes1}
Take $\lambda\in P$ and let $a\geq 2$, $b\geq 1$, $d\geq1$, and $m\geq1$.\\
$(i)$
If there exists a neighborhood of type $(a+d,b)$ in $\lambda$,
then there exist $d+1$ neighborhoods of type $(a,b)$ in $\lambda$.\\
$(ii)$
If there exists a neighborhood of type $(m(a-1)+1,mb)$ in $\lambda$,
then there exists a neighborhood of type $(a,b)$ in $\lambda$.
\end{lem}

\begin{proof}
Let $(i_1,\cdots,i_n)=w_\lambda\cdot(1,\cdots,n)$.

(i). \quad
Suppose that $(i_l,i_{l+a+d-1})$ is a neighborhood of type $(a+d,b)$
 in $\lambda$.
Set $I_j:=(i_j,i_{j+a-1})$ ($l\leq j \leq l+d$).
If there exists $j$ such that $I_j$ is not a neighborhood of type $(a,b)$, then
$\lambda_{i_j}-\lambda_{i_{j+a-1}}>b$ or
$\lambda_{i_j}-\lambda_{i_{j+a-1}}=b$ and $i_{j}<i_{j+a-1}$.
Thus, $\lambda_{i_l}-\lambda_{i_{l+a+d-1}}>b$ or 
$\lambda_{i_l}-\lambda_{i_{l+a+d-1}}=b$ and $i_l<i_{l+a+d-1}$.
However this is inconsistent with the fact that $(i_l,i_{l+a+d-1})$
 is a neighborhood of type $(a+d,b)$.
Hence any $I_j$ $(l\leq j \leq l+d)$ is a neighborhood of type $(a,b)$.

(ii). \quad
Suppose that $(i_l,i_{l+m(a-1)})$ is a neighborhood of type $(m(a-1)+1,mb)$.
Set $I_j:=(i_{l+(j-1)(a-1)},i_{l+j(a-1)})$ ($1\leq j \leq m$).
Assume that no pair $I_j$ is a neighborhood of type $(a,b)$.
Then for any $1\leq j \leq m$, we see
$\lambda_{i_{l+(j-1)(a-1)}}-\lambda_{i_{l+j(a-1)}}>b$ or
 $\lambda_{i_{l+(j-1)(a-1)}}-\lambda_{i_{l+j(a-1)}}=b$ and
 $i_{l+(j-1)(a-1)}<i_{l+j(a-1)}$.
Hence $\lambda_{i_l}-\lambda_{i_{l+m(a-1)}}>mb$ or
 $\lambda_{i_l}-\lambda_{i_{l+m(a-1)}}=mb$ and $i_l<i_{l+m(a-1)}$.
However this is inconsistent with the fact that $(i_l,i_{l+m(a-1)})$
 is a neighborhood of type $(m(a-1)+1,mb)$.
Therefore there exists at least one $I_j$
 which is a neighborhood of type $(a,b)$.
\end{proof}

\begin{lem}\label{lem:NbdOfDiffTypes2}
If $\lambda$ has a neighborhood of type $(m(k+1),m(r-1))$,
then $\lambda$ has a neighborhood of type $(k+1,r-1)$.

Let $d\geq1$.
If $\lambda$ has a neighborhood of type $(m(k+1)+d,m(r-1))$,
then $\lambda$ has two neighborhoods of type $(k+1,r-1)$.
\end{lem}

\begin{proof}
If $\lambda$ has a neighborhood of type $(m(k+1),m(r-1))$,
then from Lemma \ref{lem:NbdOfDiffTypes1} (i),
 $\lambda$ has a neighborhood of type $(mk+1,m(r-1))$.
Thus from Lemma \ref{lem:NbdOfDiffTypes1} (ii),
$\lambda$ has a neighborhood of type $(k+1,r-1)$.

If $\lambda$ has a neighborhood of type $(m(k+1)+d,m(r-1))$,
then from Lemma \ref{lem:NbdOfDiffTypes1} (i),
 $\lambda$ has a neighborhood of type $(m(k+1)+1,m(r-1))$.
Thus from Lemma \ref{lem:NbdOfDiffTypes1} (ii),
 $\lambda$ has a neighborhood of type $(k+2,r-1)$.
Therefore from Lemma \ref{lem:NbdOfDiffTypes1} (i),
$\lambda$ has two neighborhoods of type $(k+1,r-1)$.
\end{proof}

Let us show that any element in $B^{(k,r)}$ is intertwined with each other.
\begin{defn}\normalfont
Take $\lambda\in P$ and let $(i_1,\cdots,i_n)=w_\lambda\cdot(1,\cdots,n)$.
We call $\lambda'\in P$ {\it an enlargement of $\lambda$} if
$\rho(\lambda)=\rho(\lambda')$ and
$\lambda'_{i_{a}}-\lambda'_{i_{a+1}} > \max\{[\frac{n}{k+1}](r-1),\lambda_{i_{a}}-\lambda_{i_{a+1}}\}$ for any $1\leq a\leq n-1$.
\end{defn}
Note that there always exists an enlargement $\lambda'$ of $\lambda$.
We see that $\lambda'\in B^{(k,r)}$ and $\sharp^{(k,r)}(\lambda')=0$.

Let us make $\lambda'$ from $\lambda$ applying some $s_i$ and $\omega$.

\begin{prop}[enlarging procedure]
Take $\lambda\in P$ and its enlargement $\lambda'$.
Then there exists a finite sequence
$\lambda=\nu^{(0)},\nu^{(1)},\cdots,\nu^{(N-1)},\nu^{(N)}=\lambda'$
 and $(i)$, $(ii)$, or $(iii)$ holds for any $0\leq a \leq N-1$$:$\\
\quad$(i)$ $\nu^{(a+1)}=\omega^{\pm1}\nu^{(a)}$,\\
\quad$(ii)$ $\nu^{(a+1)}=s_i\nu^{(a)}$ and $\nu^{(a+1)}\succ\nu^{(a)}$,\\
\quad$(iii)$ $\nu^{(a+1)}=s_i\nu^{(a)}$, $\nu^{(a+1)}\prec\nu^{(a)}$,
 and $\nu^{(a+1)}$ is intertwined with $\nu^{(a)}$.\\
We call this sequence the enlarging procedure.
\end{prop}
\begin{proof}
Let $(i_1,\cdots,i_n)=w_\lambda\cdot(1,\cdots,n)$.
Set $M_i=\lambda'_i-\lambda_i$.
Then $M_{i_a}-M_{i_n}\geq 0$ for any $1\leq a\leq n$.
Note that $\lambda_{i_1}$ is the leftmost maximum component of $\lambda$.

Let $\Delta_{i_1}=s_{i_1}\cdots s_{n-1}\omega s_{1}\cdots s_{i_1-1}$.
Then $\Delta_{i_1}\lambda=(\lambda_1,\cdots,\lambda_{i_1}+1,\cdots,\lambda_n)$
and
 $\Delta_{i_1}^{M_{i_1}-M_{i_n}}\lambda=(\lambda_1,\cdots,\lambda_{i_1}+M_{i_1}-M_{i_n},\cdots,\lambda_n)$.
Note that for each step in
$\lambda,s_{i_1-1}\lambda,\cdots,\Delta_{i_1}\lambda,\cdots,\Delta_{i_1}^{M_{i_1}-M_{i_n}}\lambda$,
 the case (i) or (ii) holds.

Next, set $\lambda^{(1)}=\Delta_{i_1}^{M_{i_1}-M_{i_n}}\lambda$
 and apply $\Delta_{i_2}^{M_{i_2}-M_{i_n}}$ to $\lambda^{(1)}$,
Then we obtain $\lambda^{(2)}=(\lambda^{(1)}_1,\cdots,\lambda^{(1)}_{i_2}+M_{i_2}-M_{i_n},\cdots,\lambda^{(1)}_n)$.
Suppose that $s_i\nu$ and $\nu$ are serial elements satisfying $s_i\nu\prec\nu$
 in the sequence $\lambda^{(1)},\cdots,\lambda^{(2)}$.
Then we see that $s_i\nu$ and $\nu$ are intertwined.
In fact, by the definition of the procedure,
 we see $\nu_{i}-\nu_{i+1}>[\frac{n}{k+1}](r-1)$.
Hence we have $u_{\nu}(x_{i}/x_{i+1}) \neq 1, t^{\pm1}$.
Therefore in such a case, (iii) holds.

Inductively, we obtain $\lambda^{(n-1)}=(\lambda'_1-M_{i_n},\cdots,\lambda'_n-M_{i_n})$.
Similarly, if $s_i\nu$ and $\nu$ are serial elements satisfying $s_i\nu\prec\nu$
 in the sequence $\lambda,\cdots,\lambda^{(n-1)}$,
 then $s_i\nu$ and $\nu$ are intertwined, and (iii) holds.

Finally, shift it by $\omega^{nM_{i_n}}$ and
 we obtain $\lambda'$.
\end{proof}

\begin{lem}\label{lem:enlarge}
Take $\lambda\in B^{(k,r)}$ and its enlargement $\lambda'$.
Then $\lambda'$ is intertwined with $\lambda$.
\end{lem}

\begin{proof}
Consider the enlarging procedure
 $\lambda,\cdots,\nu^{(a)},\nu^{(a+1)}\cdots,\lambda'$
 given in previous proposition.

Let us show that $\nu^{(a+1)}$ is intertwined with $\nu^{(a)}$ by induction.
In the case (i) and (iii), they are intertwined.
Let us check for the case (ii).
Set $\nu=\nu^{(a)}$ and suppose $s_i\nu\succ\nu$.
By the hypothesis of induction, $\nu$ is intertwined with $\lambda$.
Thus we see $\sharp^{(k,r)}(\nu)=0$.
Assume that $u_\nu(x_{i}/x_{i+1}) =t^d$ for some $d=-1,0,1$.
Then
$(\rho(\nu)_{i+1}-\rho(\nu)_{i},\nu_{i+1}-\nu_{i})
=(m(k+1)+d,m(r-1))$ for some $m\geq 1$.
Thus $\nu$ has a neighborhood
 $(i+1,i)$ of type $(m(k+1)+d+1,m(r-1))$.
Hence from Lemma \ref{lem:NbdOfDiffTypes2},
 $\nu$ has a neighborhood of type $(k+1,r-1)$.
However it is inconsistent with $\sharp^{(k,r)}(\nu)=0$.
Therefore $s_i\nu$ and $\nu$ are intertwined.
\end{proof}

Using this lemma, we have the following statement.

\begin{lem}\label{lem:Intertwined}
Any element in $B^{(k,r)}$ is intertwined with each other.
\end{lem}
\begin{proof}
Take $\lambda,\mu\in B^{(k,r)}$, and
let $(i_1,\cdots,i_n)=w_\lambda\cdot(1,\cdots,n)$ and
$(j_1,\cdots,j_n)=w_\mu\cdot(1,\cdots,n)$.
Take $M>[\frac{n}{k+1}](r-1)+\sum_{a<b}(|\lambda_a-\lambda_b|+|\mu_a-\mu_b|)$.
Then from Lemma \ref{lem:enlarge},
 we obtain enlargements $\lambda'$ and $\mu'$ such that
 $\lambda'_{i_a}=\mu'_{j_a}=nM-aM$, $\lambda'$ is intertwined with $\lambda$,
 and $\mu'$ is intertwined with $\mu$.
Note that $\lambda'$ is a permutation of $\mu'$.
Since $|\lambda'_i-\lambda'_j|$ and $|\mu'_i-\mu'_j|>[\frac{n}{k+1}](r-1)$
 for any $1\leq i\neq j\leq n$,
 we see that $\lambda'$ and $\mu'$ are intertwined.
Therefore $\lambda$ and $\mu$ are intertwined.
\end{proof}

The well-definedness of $E_\lambda$ at (\ref{eq:Param})
 is shown by checking $Y$-eigenvalues:
\begin{lem}\label{lem:well-def}
Let $\lambda\in P$.
If there does not exist $\mu\neq\lambda$
 such that $t^{\rho(\lambda)}q^{\lambda}=t^{\rho(\mu)}q^{\mu}$
 at $(\ref{eq:Param})$,
 then $E_\lambda$ has no pole at $(\ref{eq:Param})$.
\end{lem}
\begin{proof}
Set $Y(w)=\sum_{i=1}^{n} Y_i w^{i-1}$ and
 $\lambda(w)=\sum_{i=1}^{n} t^{\rho(\lambda)_i}q^{\lambda_i} w^{i-1}$.
We see that
\bea
E_\lambda=\prod_{\mu\prec\lambda}\frac{Y(w)-\mu(w)}{\lambda(w)-\mu(w)}x^\lambda.
\eea
By the hypothesis, the right hand side does not have a pole at (\ref{eq:Param})
 for generic $w$.
Hence $E_\lambda$ has no pole at (\ref{eq:Param}).
\end{proof}

\begin{lem}\label{lem:DimOfEigen}
Let $\sharp^{(k,r)}(\lambda)\leq 1$.
Then there is no $\mu\in P$
 such that $\mu\neq\lambda$ and $t^{\rho(\mu)}q^\mu=t^{\rho(\lambda)}q^\lambda$.
\end{lem}
\begin{proof}
Suppose that there exists such an element $\mu\in P$.
Let $j_0$ be such that
$\rho(\lambda)_{j_0}=\max\{ \rho(\lambda)_i\ ;\ \rho(\lambda)_i\neq\rho(\mu)_i \}$.

Let $m_0$, $j_1$ and $m_1$ be such that
 $\rho(\mu)_{j_0}=\rho(\lambda)_{j_0}-m_0(k+1)$,
 $\rho(\mu)_{j_1}=\rho(\mu)_{j_0}+(k+1)$, and
 $\rho(\mu)_{j_1}=\rho(\lambda)_{j_1}-m_1(k+1)$.
Note that $m_0\geq1$ because of the definition of $j_0$.
Set $m_2=\frac{\rho(\lambda)_{j_0}-\rho(\lambda)_{j_1}}{k+1}$.
Note that $m_2=m_0-m_1-1\geq 1$.
In fact, if $m_1\neq 0$, then by the definition of $j_0$,
 we see that $m_2\geq1$.
If $m_1=0$, then $m_2\geq0$. Since $j_0\neq j_1$, we see $m_2\neq0$.
Hence $m_2\geq1$.

Then we have $\lambda_{j_0}-\lambda_{j_1}\geq m_2(r-1)$.
In fact, if $\lambda_{j_0}-\lambda_{j_1}< m_2(r-1)$,
 then $(j_0,j_1)$ is a neighborhood
 of type $(m_2(k+1)+1,m_2(r-1))$ in $\lambda$.
Thus from Lemma \ref{lem:NbdOfDiffTypes2}, we have $\sharp^{(k,r)}(\lambda)\geq 2$
 and this is inconsistent with $\sharp^{(k,r)}(\lambda)\leq 1$.

Since $t^{\rho(\mu)}q^\mu=t^{\rho(\lambda)}q^\lambda$, we see that
\bea
\mu_{j_1}-\mu_{j_0}&=&(\lambda_{j_1}-m_1(r-1))-(\lambda_{j_0}-m_0(r-1)) \\
&=&\lambda_{j_1}-\lambda_{j_0}+(-m_1+m_0)(r-1) \\
&\leq&(-m_2-m_1+m_0)(r-1) \\
&=&r-1.
\eea
Here, if the equality $\mu_{j_1}-\mu_{j_0} = r-1$ holds,
 then $\lambda_{j_0}-\lambda_{j_1}=m_2(r-1)$.
Hence we see $j_0<j_1$.
In fact, if $j_0>j_1$,
 then $(j_0,j_1)$ is a neighborhood of type $(m_2(k+1)+1,m_2(r-1))$ in $\lambda$.
Thus from Lemma \ref{lem:NbdOfDiffTypes2}, $\sharp^{(k,r)}(\lambda)\geq 2$,
 and this is inconsistent with $\sharp^{(k,r)}(\lambda)\leq 1$.

Hence $\mu$ has a neighborhood $(j_1,j_0)$ of type $(k+2,r-1)$.
From Lemma \ref{lem:NbdOfDiffTypes2},
 $\sharp^{(k,r)}(\mu)\geq2$.
Because $t^{\rho(\mu)}q^\mu=t^{\rho(\lambda)}q^\lambda$,
 we see that $\lambda$ also has two wheels.
However this is inconsistent with the hypothesis $\sharp^{(k,r)}(\lambda)\leq 1$.
\end{proof}

\begin{lem}\label{lem:ZerosOfSuffScattered}
Let $\lambda=(0,M,2M,\cdots,(n-1)M)$ with $M> 2[\frac{n}{k+1}](r-1)$.
Note that $\lambda\in B^{(k,r)}$ and $\sharp^{(k,r)}(\lambda)=0$.
Then the nonsymmetric Macdonald polynomial $E_\lambda$ has
 no pole at $(\ref{eq:Param})$ and we have
 $\zeta(u_0(E_\lambda))=\left[\frac{n}{k+1}\right]$.
\end{lem}

\begin{proof}
From Lemma \ref{lem:well-def} and \ref{lem:DimOfEigen},
we see that $E_\lambda$ has no pole at (\ref{eq:Param}).

Let us compute $\zeta(u_0(E_\lambda))$.
From Lemma \ref{lem:u_0formula_omega} and \ref{lem:u_0formula_s_i},
 we obtain $u_0(E_\lambda)$ from $u_0(E_{0})=1$.
Let us find factors of the form $(t^{m(k+1)}q^{m(r-1)}-1)$.
We have
\bea
&&u_0(E_{(0,0,\cdots,0,0,(n-1)M)})=u_0(E_{(0,0,\cdots,0,0,0)})
c_1\prod_{i=1}^{M-1}\prod_{j=n-1}^{1}\frac{t^{j+1}q^i-1}{t^jq^i-1}, \\
&&u_0(E_{(0,0,\cdots,0,(n-1)M,nM)})=u_0(E_{(0,0,\cdots,0,0,(n-1)M)})
c_2\prod_{i=1}^{M-1}\prod_{j=n-2}^{1}\frac{t^{j+1}q^i-1}{t^jq^i-1}, \\
&&\cdots, \\
&&u_0(E_{(0,M,2M,\cdots,(n-1)M)})=u_0(E_{(0,0,2M,\cdots,(n-1)M)})
c_{n-1}\prod_{i=1}^{M-1}\prod_{j=1}^{1}\frac{t^{j+1}q^i-1}{t^jq^i-1}.
\eea
Here, $c_i$ is a factor which does not contain $(t^iq^N-1)$ with $N< M$.
Hence 
\bea
u_0(E_\lambda)&=&u_0(E_{(0,M,2M,\cdots,(n-1)M)}) \\
&=&c\prod_{l=1}^{n-1}\prod_{i=1}^{M-1}\prod_{j=1}^{n-l}\frac{t^{j+1}q^i-1}{t^jq^i-1} \quad\mbox{where $c=\prod c_i$},\\
&=&c\prod_{l=1}^{n-1}\prod_{i=1}^{M-1}\frac{t^{n-l+1}q^i-1}{tq^i-1} \\
&=&c\prod_{l=2}^{n}\prod_{i=1}^{M-1}\frac{t^{l}q^i-1}{tq^i-1}.
\eea
Therefore, the multiplicity of the factor $(t^{k+1}q^{r-1}-1)$
 is $[\frac{n}{k+1}]$.
\end{proof}

\begin{lem}\label{lem:ZerosOfNoWheel}
For any $\lambda\in B^{(k,r)}$,
the nonsymmetric Macdonald polynomial $E_\lambda$ has no pole at $(\ref{eq:Param})$
 and we have
$\zeta(u_0(E_{\lambda}))=[\frac{n}{k+1}]$.
\end{lem}
\begin{proof}
This is clear from Lemma \ref{lem:Intertwined} and \ref{lem:ZerosOfSuffScattered}.
\end{proof}

\begin{lem}\label{lem:ZerosOfOneWheel}
For $\mu\in P$,
let $(i_1,\cdots,i_n)=w_\mu\cdot(1,\cdots,n)$.
Suppose that $\mu$ satisfies the following condition$:$
\bean
&&\begin{array}{l}
\mbox{for some $1\leq l \leq n-k$,}\\
\mbox{$\mu$ has a neighborhood $(i_l,i_{l+k})$ of type $(k+1,r-1)$, and} \\
\mbox{$|\mu_i-\mu_j|> 2[\frac{n}{k+1}](r-1)$}\\
\qquad\qquad\mbox{for any $1\leq i\neq j\leq n$ except for $(i,j)=(i_a,i_b)$ $(l\leq  a,b\leq l+k)$}.
\end{array}\label{eq:mu_onewheel}
\eean
Note that $\sharp^{(k,r)}(\mu)=1$.

Then, $E_\mu$ has no pole at $(\ref{eq:Param})$
 and $\zeta(u_0(E_\mu))=[\frac{n}{k+1}]-1$.
\end{lem}
\begin{proof}
From Lemma \ref{lem:well-def} and \ref{lem:DimOfEigen},
 we see that $E_\mu$ has no pole at (\ref{eq:Param}).

Let us compute $\zeta(u_0(E_\mu))$.
Fix an enlargement $\mu'$ of $\mu$, and consider the enlarging procedure.
Note that $\sharp^{(k,r)}(\mu)=1$ and $\sharp^{(k,r)}(\mu')=0$.
Thus $\mu$ and $\mu'$ are not intertwined.

Suppose that $\nu$ and $s_i\nu$
 are serial elements in the enlarging procedure such that
 $\mu$ and $\nu$ are intertwined and $\nu$ is not intertwined with $s_i\nu$.
Then we see that
 $\sharp^{(k,r)}(\nu)=1$,
 $u_{\nu}(x_i/x_{i+1})=t^d$ for some $d=-1,0,1$,
 $s_i\nu\succ\nu$,
 and $\nu$ satisfies the condition (\ref{eq:mu_onewheel}) replaced by $\nu$.

Let us show that $u_{\nu}(x_i/x_{i+1})=t^{-1}$.
Assume that $u_{\nu}(x_i/x_{i+1})=t^d$ for $d=0,1$.
Then $(\rho(\nu)_{i+1}-\rho(\nu)_{i},\nu_{i+1}-\nu_{i})=(m(k+1)+d,m(r-1))$
 for some $m\geq 1$.
Thus $\nu$ has a neighborhood $(i+1,i)$ of type $(m(k+1)+d+1,m(r-1))$,
and from Lemma \ref{lem:NbdOfDiffTypes2},
 $\nu$ has two neighborhoods of type $(k+1,r-1)$.
However this is inconsistent with $\sharp^{(k,r)}(\nu)=1$.

Let us show that $\sharp^{(k,r)}(s_i\nu)<\sharp^{(k,r)}(\nu)$.
We have shown that $u_{\nu}(x_i/x_{i+1})=t^{-1}$.
Then $(\rho(\nu)_{i+1}-\rho(\nu)_i,\nu_{i+1}-\nu_i)=(m(k+1)-1,m(r-1))$
 for some $m\geq 1$.
Here we see $m=1$.
In fact, if $m\geq2$, then
 from the relation $\nu_{i+1}-\nu_i=m(r-1)$ and (\ref{eq:mu_onewheel}),
 we see $m\geq 2[\frac{n}{k+1}]$.
Hence $\rho(\nu)_{i+1}-\rho(\nu)_i\geq n$.
However this is contradiction.
Therefore $\nu$ has a neighborhood $(i+1,i)$ of type $(k+1,r-1)$
 such that $\nu_{i+1}-\nu_{i}=r-1$.
Let $(j_1,\cdots,j_n)=w_\nu\cdot(1,\cdots,n)$.
Then $j_l=i+1$ and $j_{l+k}=i$ for some $1\leq l\leq n-k$.
We see that for some $l\leq a\leq l+k$,
 $J=(j_{a},\cdots,j_{l+k},j_l,\cdots,j_{a-1})$ is a wheel in $\nu$.
Moreover, by the definition of wheels, $s_iJ$ is not a wheel in $s_i\nu$.
Therefore from Remark \ref{rem:GenAndBrakeOfWheel},
 we have $\sharp^{(k,r)}(s_i\nu)<\sharp^{(k,r)}(\nu)$.

Since $\sharp^{(k,r)}(\nu)=1$, we obtain $\sharp^{(k,r)}(s_i\nu)=0$.
Let us show that $s_i\nu$ and $\mu'$ are intertwined.
Let $\xi$ and $s_j\xi$ are any serial elements between $s_i\nu$ and $\mu'$
 in the enlarging procedure.
If $s_j\xi\prec\xi$, then $\xi$ and $s_j\xi$ are intertwined.
Suppose that $s_j\xi\succ\xi$ and $\sharp^{(k,r)}(\xi)=0$.
Assume that $u_{\xi}(x_j/x_{j+1})=t^d$ for some $d=-1,0,1$.
Then $(\rho(\xi)_{j+1}-\rho(\xi)_{j},\xi_{j+1}-\xi_{j})=(m(k+1)+d,m(r-1))$
 for some $m\geq 1$.
Thus $\xi$ has a neighborhood $(j+1,j)$ of type $(m(k+1)+d+1,m(r-1))$,
and from Lemma \ref{lem:NbdOfDiffTypes2},
 $\xi$ has a neighborhood of type $(k+1,r-1)$.
However this is inconsistent with $\sharp^{(k,r)}(\xi)=0$.
Hence $\xi$ and $s_j\xi$ are intertwined and $\sharp^{(k,r)}(s_j\xi)=0$.
Inductively, we obtain that $s_i\nu$ and $\mu'$ are intertwined.

Consequently, we obtain
\bea
\zeta(u_0(E_{\mu}))
&=&\zeta(u_0(E_{\nu}))\\
&=&\zeta(u_0(E_{s_i\nu}))-1\\
&=&\zeta(u_0(E_{\mu'}))-1 \\
&=&\left[\frac{n}{k+1}\right]-1.
\eea
\end{proof}

Now we give a lower estimate of $I^{(k,r)}$.

\begin{prop}\label{prop:LowerEstimate}
For any $\lambda\in B^{(k,r)}$,
 the nonsymmetric Macdonald polynomial $E_\lambda$
 has no pole at $(\ref{eq:Param})$.
Moreover, $E_\lambda$ specialized at $(\ref{eq:Param})$ belongs to $I^{(k,r)}$.
\end{prop}

\begin{proof}
Let $\lambda\in B^{(k,r)}$.
Then from Lemma \ref{lem:ZerosOfNoWheel},
the well-definedness of $E_\lambda$ is proved
and we have $\zeta(u_0(E_\lambda))=\left[\frac{n}{k+1}\right]$.

Let us show that $E_\lambda\in I^{(k,r)}$.
Take $D > 2(\max_i\{ |\lambda_i| \}+1)$, $M> 2[\frac{n}{k+1}](r-1)$,
 and $N>2M+2D$.
Take $\mu^{(0)}\in S'{}^{(k,r)}$
 satisfying the following condition:
\bean
&&\begin{array}{l}
\mbox{for some $1\leq l \leq n-k$,}\\
\mbox{$\mu^{(0)}$ has a neighborhood $(i_l,i_{l+k})$ of type $(k+1,r-1)$}\\
\mbox{such that $\mu^{(0)}_{i_l}-\mu^{(0)}_{i_{l+k}}\leq r-2$,
 and $|\mu^{(0)}_i-\mu^{(0)}_j|> M$}\\
\mbox{ for any $1\leq i\neq j\leq n$ except for $(i,j)=(i_a,i_b)$ ($l\leq  a,b\leq l+k$)}.
\end{array}\label{eq:mu^0}
\eean
Here we set $(i_1,\cdots,i_n)=w_{\mu^{(0)}}\cdot(1,\cdots,n)$.
We see $\sharp^{(k,r)}(\mu^{(0)})=1$.
Define a finite set
\bea
S(\mu^{(0)})=\{ \mu\in P ;
 \mbox{$0\leq \mu_i-\mu^{(0)}_i \leq D$ for $i\neq i_l,i_{l+1},\cdots,i_{l+k}$} \}.
\eea
Then from Lemma \ref{lem:ZerosOfOneWheel}, for any $\mu\in S(\mu^{(0)})$,
we see that $E_\mu$ has no pole at (\ref{eq:Param})
 and $\zeta(u_0(E_\mu))=[\frac{n}{k+1}]-1$.
From the duality relation, we have
\[
u_\mu(E_\lambda)=\frac{u_\lambda(E_\mu)}{u_0(E_\mu)}u_0(E_\lambda).
\]
Hence we see that $\zeta(u_\mu(E_\lambda))\geq 1$.
Therefore, from the same argument as Lemma \ref{lem:EquivDef},
we have that $E_\lambda$ specialized at (\ref{eq:Param}) belongs to $I^{(k,r)}$.
\end{proof}

\section{Irreducibility of $I^{(k,r)}$}

In this section, we give an upper estimate of the character of $I^{(k,r)}$,
 and thereby, we show that $I^{(k,r)}$ is irreducible.
Thus, we will complete the proof of Theorem \ref{thm:NoWheel}.

First, we give an upper estimate of the character of $I^{(k,r)}$.
Recall the definition of $I^{(k,r)}$:
\[
I^{(k,r)}=\{f\in V ;\ f(z)=0\ \mbox{for any $z\in Z^{(k,r)}$} \}
\]
where
\bea
Z^{(k,r)}=\{ z\in\mathbb{K}^n&;&
\exists i_1,\cdots ,i_{k+1},
\exists s_1,\cdots,s_{k}\in\mathbb{Z}_{\geq0} \\
&&\mbox{ such that $z_{i_{a+1}}=z_{i_a}tq^{s_a}$ for $1\leq a\leq k$}, \\
&& \quad\mbox{ $\sum_a s_a\leq r-2$, and $i_{a}<i_{a+1}$ if $s_a=0$}\}.
\eea
We define the character of $I^{(k,r)}$ as follows.
For any Laurent polynomial $f=\sum_\lambda c_\lambda x^\lambda$,
 we set $\deg(f)=\max\{|\lambda_i| ; c_\lambda\neq0, 1\leq i \leq n \}$.
We introduce the filtration
\bea
V_{(d)}&=&\{f\in V; \deg(f)\leq d\},
\eea
and define $I^{(k,r)}_{(d)}= I^{(k,r)} \cap V_{(d)}$.

Let us take the limit $u\rightarrow 1$.
Namely, take $V_{\mathrm{reg}}=\{f\in V; \mbox{$f$ is regular at $u=1$} \}$ and
let $V_0$ be the image of $V_{\mathrm{reg}}$ by the specialization $u=1$.
Note that $V_0=\mathbb{C}[x_1^{\pm1},\cdots,x_n^{\pm1}]$.
Then the ideal reduces to
\[
{I^{(k,r)}}'=\{f\in V_0 ;\ f(z)=0 \ \mbox{for any $z\in {Z^{(k,r)}}'$} \}
\]
where
\bea
{Z^{(k,r)}}'=\{ z\in\mathbb{C}^n&;&
\exists j_1<\cdots <j_{k+1},
\exists p_1,\cdots p_{k+1}\in\mathbb{Z}_{\geq0},
\exists w\in\mathbb{C} \\
&&\mbox{ such that } z_{j_a}=\tau^{p_a}w \mbox{ for } 1\leq a\leq k+1 \\
&&\quad\mbox{ and } 0\leq p_a\leq r-2 \}.
\eea
The correspondence between $Z^{(k,r)}$ and ${Z^{(k,r)}}'$ is as follows.
For $z\in Z^{(k,r)}$, we see that
$z_{i_a}=z_{i_1}t^{a-1}q^{s_1+\cdots+s_{a-1}}$.
Thus for $\sigma\in\mathfrak{S}_{k+1}$ such that
 $i_{\sigma(1)}<\cdots<i_{\sigma(k+1)}$,
 we take $j_a=i_{\sigma(a)}$, $p_a=p_{\sigma^{-1}(1)}+s_1+\cdots+s_{\sigma(a)-1}$,
 $p_{\sigma^{-1}(1)}$ and $w$ satisfying $\tau^{p_{\sigma^{-1}(1)}}w=z_{i_1}$.
We also introduce the filtration
\bea
(V_0)_{(d)}&=&\{f\in V_0; \deg(f)\leq d\},
\eea
and set ${I^{(k,r)}_{(d)}}'={I^{(k,r)}}' \cap (V_0)_{(d)}$.
Note that the character of the ideal does not decrease under this limit.
Namely $\dim_{\mathbb{K}} I^{(k,r)}_{(d)}\leq \dim_{\mathbb{C}} {I^{(k,r)}_{(d)}}'$.

Fix an arbitrary non-negative integer $M\geq 0$ and
let us estimate the dimension of ${I^{(k,r)}_{(M)}}'$.
We set $I_M={I^{(k,r)}_{(M)}}'$,
 $P_M=\{ \lambda\in P ; -M\leq\lambda_i\leq M \}$,
 $S^{(k,r)}_M= S^{(k,r)} \cap P_M$,
 and $B^{(k,r)}_M= B^{(k,r)} \cap P_M$.
Consider the tensor algebra $T(\mathrm{span}_{\mathbb{C}}\{e_d; -M\leq d \leq M\})$
and denote it by $R_M$.
Denote its $n$-th tensor subspace by $R_{M,n}$.
For simplicity, we write $e_{\lambda_1}\cdots e_{\lambda_n}
=e_{\lambda_1}\otimes \cdots\otimes e_{\lambda_n}$ and
 $e_{\lambda}=e_{\lambda_1}\cdots e_{\lambda_n}$.
Define
\[ e(w)_M=\sum_{-M\leq d \leq M} e_dw^d, \]
then
\[
e(z_1)_M\cdots e(z_n)_M
=\sum_{\lambda\in P_M} e_\lambda z^\lambda.
\]
We introduce a non-degenerate pairing
\bea
\langle\cdot,\cdot\rangle:R_{M,n}\times (V_0)_{(M)}
 \rightarrow \mathbb{C}
\eea
 by $\langle e_\lambda , x^\mu \rangle=\delta_{\lambda,\mu}$.
Then \[ \langle e(z_1)_M\cdots e(z_n)_M,f \rangle =f(z_1,\cdots,z_n) \]
for any $f\in (V_0)_{(M)}$.
By this pairing, the ideal $I_M$ is written as follows:
\[ R_{M,n}/J_M \cong I_M. \]
Here, 
\bea
J_M&=&\mathrm{span}_\mathbb{C}
\{  e(z_1)_M \stackrel{\scriptsize i_1 th}{\cdots e(\tau^{p_1}w)_M\cdots}
 \stackrel{i_{k+1} th}{e(\tau^{p_{k+1}}w)_M}\cdots e(z_n)_M ;\\
&&\qquad\qquad\qquad\qquad
 i_1<\cdots<i_{k+1},0\leq p_i\leq r-2, w,z_i\in\mathbb{C} \}. 
\eea

Let us give a spanning set of the quotient space $R_{M,n}/J_M$.

\begin{prop}
A spanning set of $R_{M,n}/J_M$ is given by the quotient image of the set
 $\{e_\lambda ; \lambda\in B^{(k,r)}_M \}$.
In other words, in $R_{M,n}/J_M$,
 the image of $e_\lambda$ $(\lambda\in S^{(k,r)}_M)$ is
written as a linear combination of the image of $e_\mu$
 $(\mu\in B^{(k,r)}_M)$.
\end{prop}

\begin{proof}
First we calculate elements of $J_M$.
We define some notations.
\[
\mathbb{Z}^{k+1}_{r-1}=\{(\eta_1,\cdots,\eta_{k+1})\in\mathbb{Z}^{k+1};
 0\leq \eta_i\leq r-2 \mbox{ for $1\leq i \leq k+1$} \}.
\]
For $\eta\in\mathbb{Z}^{k+1}_{r-1}$ and an integer $d\in\mathbb{Z}$, define
\[
S(\eta,d)=\{\eta'\in\mathbb{Z}^{k+1};
 \sum_i \eta_i=d, \eta'_i=\eta_i \mbox{ mod } r-1\} .
\]

Fix $1\leq j_1<\cdots<j_{k+1}\leq n$.
Then,
\bea
&&e(z_1) \stackrel{j_1 th}{\cdots e(\tau^{p_1}w)\cdots}
 \stackrel{j_{k+1} th}{e(\tau^{p_{k+1}}w)}\cdots e(z_n) \\
&&\ =\sum_{\nu\in\mathbb{Z}^n}
 e_\nu \tau^{p_1\nu_{j_1}+\cdots+p_{k+1}\nu_{j_{k+1}}}
 w^{\nu_{j_1}+\cdots \nu_{j_{k+1}}}
 z_1^{\nu_1} \cdots \widehat{z_{j_1}^{\nu_{j_1}}}\cdots
 \widehat{z_{j_{k+1}}^{\nu_{j_{k+1}}}} \cdots z_n^{\nu_n}\\
&&\ =\sum_{\eta,d,\hat{\nu}} r_{(j_a),\eta,d,\hat{\nu}}
 \tau^{\sum_{a=1}^{k+1} p_a \eta_a} w^d
 z_1^{\nu_1} \cdots \widehat{z_{j_1}^{\nu_{j_1}}}\cdots
 \widehat{z_{j_{k+1}}^{\nu_{j_{k+1}}}} \cdots z_n^{\nu_n}.
\eea
Here, the last sum runs over
$\eta\in\mathbb{Z}^{k+1}_{r-1}$,
$d\in\mathbb{Z}$,
$\hat{\nu}=(\nu_1,\cdots,\widehat{\nu_{j_1}},\cdots,\widehat{\nu_{j_{k+1}}},\cdots,\nu_n)\in\mathbb{Z}^{n-k-1}$,
 and $r_{(j_a),\eta,d,\hat{\nu}}$ is given by
\[
r_{(j_a),\eta,d,\hat{\nu}}=\sum_\lambda e_\lambda.
\]
The sum runs over $\lambda\in P_M$ such that
$(\lambda_{j_1},\cdots,\lambda_{j_{k+1}})\in S(\eta,d)$
 and $\lambda_i=\nu_i$ for $i\neq j_1,\cdots,j_{k+1}$.

Hence in $R_{M,n}/J_M$,
\bea
 0=\sum_{\eta\in\mathbb{Z}^{k+1}_{r-1}} r_{(j_a),\eta,d,\hat{\nu}}
 \tau^{\sum_{a=1}^{k+1} p_a \eta_a}
\eea
for any $d\in\mathbb{Z}$ and $\hat{\nu}\in\mathbb{Z}^{n-k-1}$.
Since this equality holds for any $0\leq p_1,\cdots,p_{k+1}\leq r-2$, we have
\[ r_{(j_a),\eta,d,\hat{\nu}}=0 \]
for any $j_1<\cdots<j_{k+1}$,
$\eta\in\mathbb{Z}^{k+1}_{r-1}$,
$d\in\mathbb{Z}$, and
$\hat{\nu}\in\mathbb{Z}^{n-k-1}$.

Let us introduce the total ordering on $P$.
For $\lambda,\mu\in P$,
let $(i_1,\cdots,i_n)=w_\lambda(1,\cdots,n)$
 and $(j_1,\cdots,j_n)=w_\mu(1,\cdots,n)$.
We define $\lambda>'\mu$ if there exists $1\leq l\leq n$ such that
 $\lambda_{i_a}=\mu_{j_a}$ and $i_a=j_a$ for any $1\leq a \leq l-1$,
 and $\lambda_{i_l}>\mu_{j_l}$, or $\lambda_{i_l}=\mu_{j_l}$ and $i_l<j_l$.
We induce the ordering to monomials.
Namely, we define $e_\lambda>'e_\mu$ if $\lambda>'\mu$.

Suppose $\lambda\in S^{(k,r)}_M$.
Let us rewrite $e_\lambda$  in $R_{M,n}/J_M$
 as a linear combination of greater monomials with respect to $>'$.
Let $(i_1,\cdots,i_n)=w_\lambda\cdot(1,\cdots,n)$ and
 $(i_l,i_{l+k})$ be a neighborhood of type $(k+1,r-1)$ in $\lambda$.
Take $\sigma\in \mathfrak{S}_{k+1}$
 such that $i_{\sigma(l)}<i_{\sigma(l+1)}<\cdots <i_{\sigma(l+k)}$
 and let $\lambda'=(\lambda_{\sigma(i_l)},\lambda_{\sigma(i_{l+1})},\cdots,\lambda_{\sigma(i_{l+k})})$.
Take $\eta\in\mathbb{Z}^{k+1}_{r-1}$ satisfying $\eta_i\equiv\lambda'_i$
 mod $r-1$, and let $d=\sum_{i=1}^{k+1} \lambda'_i$.
Then for any $\eta'\in S(\eta,d)$ such that $\eta'\neq\lambda'$,
 we have $\eta'>'\lambda'$.
Hence we have
\[ 0=r_{(\sigma(i_a)),\eta,d,\hat{\lambda}}=e_\lambda+\sum_\mu e_\mu. \]
The last sum runs over $\mu\in P_M$ such that $\mu>'\lambda$,
 $(\mu_{\sigma(i_1)},\cdots,\mu_{\sigma(i_{k+1})})\in S(\eta,d)$
 and $\mu_i=\lambda_i$ for $i\neq i_l,\cdots,i_{l+k}$.

Continue this procedure.
Because the set $P_M$ is finite, the procedure stops in finite times.
Therefore for any $\lambda\in S^{(k,r)}_M$, we can rewrite $e_\lambda$
 as a linear combination of $e_\mu$ ($\mu\in B^{(k,r)}_M$).
\end{proof}

As a corollary, we obtain an upper estimate of the character of $I^{(k,r)}$.
\begin{cor}\label{cor:UpperEstimate}
For any $M\geq0$,
there exists a spanning set of $I^{(k,r)}_M$
 which is labeled by $B^{(k,r)}_M$.
\end{cor}

Since the upper and the lower estimates coincide,
 $\{E_\lambda;\lambda\in B^{(k,r)}\}$
 is a basis of $I^{(k,r)}$.
Now we are going to finish the proof of Theorem \ref{thm:NoWheel}.

\begin{proof}[End of the proof of Theorem \ref{thm:NoWheel}]
Let us show that any non-zero element $v\in I^{(k,r)}$ is cyclic.
From Lemma \ref{lem:DimOfEigen}, all $Y$-eigenvalues in $I^{(k,r)}$ are different.
Thus a certain $E_\lambda$ ($\lambda\in B^{(k,r)}$)
 is contained in $\HH^{(k,r)}v$.
Take another $\mu\in B^{(k,r)}$.
Then $\mu$ is intertwined with $\lambda$ from Lemma \ref{lem:Intertwined}.
Hence by applying intertwiners $A$ and $B_i$ on $E_\lambda$,
the vector $E_\mu$ is also contained in $\HH^{(k,r)}v$.
Since $\{E_\mu; \mu\in B^{(k,r)}\}$ is a basis of $I^{(k,r)}$,
 we obtain that $v$ is cyclic.
\end{proof}

\section{A series of subrepresentations defined by multi-wheel condition}

In this section, we construct a series of subrepresentations of $\HH^{(k,r)}$,
 in which $I^{(k,r)}$ appears as its member.

\begin{defn}\normalfont
Define
\bea
Z^{(k,r)}_m=\{ z\in\mathbb{K}^n;
&&\mbox{there exist $1\leq i_{l,a}\leq n$ and $s_{l,a}\in\mathbb{Z}_{\geq0}$}\\
&&\mbox{for $1\leq l \leq m$ and $1\leq a \leq k+1$ such that}\\
&&\mbox{ $i_{l,a}$ are distinct, $z_{i_{l,a+1}}=z_{i_{l,a}}tq^{s_{l,a}}$,}\\
&&\mbox{$\sum_{a=1}^{k} s_{l,a}\leq r-2$,
and $i_{l,a+1}<i_{l,a}$ if $s_{l,a}=0$} \}.
\eea
We define the ideal
\[
I_m^{(k,r)}=\{f\in V ; f(z)=0 \mbox{ for any $z\in Z^{(k,r)}_m $} \}.
\]
We call the defining condition of $I_m^{(k,r)}$ the {\it multi-wheel condition}.
\end{defn}

By the definition, we see that
\bea
&&I^{(k,r)}=I^{(k,r)}_1 \subset I^{(k,r)}_2 \subset I^{(k,r)}_3
 \subset\cdots\\
&&\quad\cdots\subset I^{(k,r)}_{m_0} \subset I^{(k,r)}_{m_0+1} = I^{(k,r)}_{m_0+2}
 = \cdots =V
\eea
where $m_0=[\frac{n}{k+1}]$.

Similarly to the single wheel case,
 we have an alternative definition of $I_m^{(k,r)}$.

Take $\lambda\in P$ and let $(i_1,\cdots,i_n)=w_\lambda\cdot(1,\cdots,n)$.
We say that
 {\it neighborhoods $(i_{l},i_{l+k})$ and $(i_{l'},i_{l'+k})$ are distinct}{}
 if $\{i_{l},i_{l+1},\cdots,i_{l+k}\}$ and $\{i_{l'},i_{l'+1},\cdots,i_{l'+k}\}$
 are disjoint.

\begin{prop}
Define the set $S_m^{(k,r)}=\{ \lambda\in P;$
 $\lambda$ has distinct $m$ neighborhoods of type $(k+1,r-1)$ $\}$.
Then the ideal $I_m^{(k,r)}$ coincides with
\bea
\{f\in V ; u_\lambda(f)=0 \mbox{ for any $\lambda\in S^{(k,r)}_m $} \}.
\eea
\end{prop}
\begin{proof}
It is proved in the same way as the single wheel case.
\end{proof}

We have the following theorem.
\begin{thm}\label{thm:MultiWheel}
The ideal $I^{(k,r)}_m$ is a representation of $\HH^{(k,r)}$.
\end{thm}

\begin{proof}
Similarly to the single wheel case,
it is sufficient to show that $\omega^{\pm1} f$, $T_i f$$\in I^{(k,r)}_m$.

Take an element $f\in I^{(k,r)}_m$ and an integer $N > 2(\deg(f)+1)$.
Let $M > 2N+2[\frac{n}{k+1}](r-1)$.

Take $\lambda^{(0)}\in S^{(k,r)}_m$ satisfying that
$(i_{l,1},i_{l,k+1})$ is a neighborhood of type $(k+1,r-1)$
 for $1\leq l \leq m$, and
$|\lambda^{(0)}_i-\lambda^{(0)}_j|> M$
 for any $1\leq i\neq j\leq n$ except for $(i,j)=(i_{l,a},i_{l,b})$
 ($1\leq l \leq m$, $1\leq  a,b\leq k+1$).
Note that $\sharp^{(k,r)}(\lambda^{(0)})=m$.
Define a finite set
\bea
S_m(\lambda^{(0)})=\{ \lambda\in P ; \mbox{$0\leq \lambda_i-\lambda^{(0)}_i \leq N$ for $i\neq i_{l,a}$ ($1\leq l \leq m,1\leq a\leq k+1$)} \}.
\eea
Note that $S_m(\lambda^{(0)})\subset S_m^{(k,r)}$.

Then, similarly to the single wheel case,
we can show that
 $u_\lambda(\omega^{\pm1} f)=u_\lambda(T_i f)=0$ ($1\leq i \leq n-1$)
 for any $\lambda\in S_m(\lambda^{(0)})$,
and we see $\omega^{\pm1} f$, $T_i f$$\in I^{(k,r)}_m$.
\end{proof}

\begin{conj}
The quotient representations $I^{(k,r)}_m/I^{(k,r)}_{m-1}$ are irreducible.
\end{conj}

We have the following statement in the case $n=k+1$.
\begin{thm}
The quotient representation $V/I_1^{(n-1,r)}$ is irreducible.
For any $\lambda \in S^{(n-1,r)}$,
the nonsymmetric Macdonald polynomial $E_\lambda$ has no pole at $(\ref{eq:Param})$.
A basis of $V/I_1^{(n-1,r)}$ is given by
 $\{ E_\lambda; \lambda \in S^{(n-1,r)} \}$ specialized at $(\ref{eq:Param})$.
\end{thm}
\begin{proof}
In the case $n=k+1$, by the definition of wheels,
 we see that $\sharp^{(n-1,r)}(\lambda)\leq 1$ for any $\lambda\in P$.
Hence from Lemma \ref{lem:well-def} and \ref{lem:DimOfEigen},
 $E_\lambda$ has no pole at $(\ref{eq:Param})$
 for any $\lambda \in S^{(n-1,r)}$.

Since the basis of $I_1^{(n-1,r)}$ is given by
 $\{ E_\lambda; \lambda \in B^{(n-1,r)} \}$
 specialized at $(\ref{eq:Param})$,
we see that the basis of $V/I_1^{(n-1,r)}$ is given by
 $\{ E_\lambda; \lambda \in S^{(n-1,r)} \}$.

Let us show that
 $(0,\cdots,0)\in S^{(n-1,r)}$ is intertwined with any $\lambda \in S^{(n-1,r)}$.
Fix $\lambda\in S^{(n-1,r)}$
 and let $(i_1,\cdots,i_n)=w_\lambda\cdot(1,\cdots,n)$.
Then $\lambda_{i_1}-\lambda_{i_n}\leq r-2$,
 or $\lambda_{i_1}-\lambda_{i_n}=r-2$ and $i_n<i_1$.

Suppose that $\lambda_{i_1}-\lambda_{i_n}\leq r-2$.
Set $\Delta_{i}=s_{i}\cdots s_{n-1}\omega s_{1}\cdots s_{i-1}$.
Then $\lambda$ is obtained by
\bea
\omega^{n\lambda_{i_n}}\Delta_{i_{n-1}}^{(\lambda_{i_{n-1}}-\lambda_{i_n})}\cdots
\Delta_{i_1}^{(\lambda_{i_1}-\lambda_{i_n})}\cdot(0,\cdots,0).
\eea
Note that for any serial elements $\nu$ and $s_i\nu$
 in the sequence $(0,\cdots,0),\cdots,\lambda$,
 we see that $\max_{j \neq j'}\{|\nu_j-\nu_j'|\}\leq r-2$.
Hence $u_\nu(x_i/x_{i+1})\neq 1,t^{\pm1}$, and
 $\nu$ is intertwined with $s_i\nu$.
Therefore $\lambda$ is intertwined with $(0,\cdots,0)$.

Suppose that $\lambda_{i_1}-\lambda_{i_n}=r-1$ and $i_n<i_1$.
Then by applying $\omega$ for some times, we obtain $\lambda'$
 such that $\lambda'_{j_1}-\lambda'_{j_n}\leq r-2$.
Here we use $(j_1,\cdots,j_n)=w_{\lambda'}\cdot(1,\cdots,n)$.
Hence $\lambda$ is intertwined with $(0,\cdots,0)$.

Therefore, from the definition of "intertwined",
 we see that $E_\lambda$ is a cyclic vector of $V/I_1^{(n-1,r)}$
 for any $\lambda \in S^{(n-1,r)}$.
\end{proof}

Combining this theorem and Theorem \ref{thm:NoWheel},
 we see that the conjecture is true in the case $n=k+1$.

Although we do not give a proof here,
 we can show that $I^{(1,r)}_2/I^{(1,r)}_1$ is irreducible
 and it is not $Y$-semisimple.
We can also construct an explicit basis of $I^{(1,r)}_2$ in terms of
 a linear combination of nonsymmetric Macdonald polynomials
 specialized at (\ref{eq:Param}).


\end{document}